\documentclass[11pt]{article}

\usepackage{amsmath,amsfonts,amsthm,amssymb,graphicx,tikz,tikz-cd,url,hyperref}
\usepackage[all,2cell,ps]{xy}

\tikzset{
math to/.tip={Glyph[glyph math command=rightarrow]},
loop/.tip={Glyph[glyph math command=looparrowleft, swap]},
loop'/.tip={Glyph[glyph math command=looparrowleft]},
 weird/.tip={Glyph[glyph math command=Rrightarrow, glyph length=1.5ex]},
  pi/.tip={Glyph[glyph math command=pi, glyph length=1.5ex, glyph axis=0pt]},
}

\theoremstyle{plain}
\newtheorem{thm}{Theorem}[section]

\newtheorem{lem}[thm]{Lemma}
\newtheorem{prop}[thm]{Proposition}

\newtheorem{cor}[thm]{Corollary}

\theoremstyle{definition}

\newtheorem{rem}[thm]{Remark}

\theoremstyle{remark}

\newcommand{\bbB}{\mathbb{B}}
\newcommand{\bbC}{\mathbb{C}}

\newcommand{\bbH}{\mathbb{H}}

\newcommand{\bbN}{\mathbb{N}}

\newcommand{\bbP}{\mathbb{P}}
\newcommand{\bbQ}{\mathbb{Q}}
\newcommand{\bbR}{\mathbb{R}}

\newcommand{\bbZ}{\mathbb{Z}}

\newcommand{\calO}{\mathcal{O}}

\newcommand{\frakA}{\mathfrak{A}}

\newcommand{\frakn}{\mathfrak{n}}

\newcommand{\al}{\alpha}
\newcommand{\gam}{\gamma}
\newcommand{\Gam}{\Gamma}

\newcommand{\Del}{\Delta}

\newcommand{\Lam}{\Lambda}
\newcommand{\sig}{\sigma}

\newcommand{\Om}{\Omega}

\DeclareMathOperator{\PSL}{PSL}

\DeclareMathOperator{\Id}{Id}

\newcommand{\bs}{\backslash}

\newcommand{\lra}{\longrightarrow}

\newcommand{\ssm}{\smallsetminus}

\newcommand{\wh}{\widehat}
\newcommand{\wt}{\widetilde}

\newenvironment{pf}{\begin{proof}}{\end{proof}}

\usepackage{tikz-cd}
\allowdisplaybreaks

\makeatletter
\let\@@pmod\pmod
\DeclareRobustCommand{\pmod}{\@ifstar\@pmods\@@pmod}
\def\@pmods#1{\mkern4mu({\operator@font mod}\mkern 6mu#1)}
\makeatother

\usetikzlibrary{arrows}

\title{One-cusped complex hyperbolic $2$-manifolds}
\author{Martin Deraux \and Matthew Stover}
\date{Dec.\ 4, 2025}

\begin{document}

\maketitle

\begin{abstract}
This paper builds one-cusped complex hyperbolic $2$-manifolds by an explicit geometric construction. Specifically, for each odd $d \ge 1$ there is a smooth projective surface $Z_d$ with $c_1^2(Z_d) = c_2(Z_d) = 6d$ and a smooth irreducible curve $E_d$ on $Z_d$ of genus one so that $Z_d \ssm E_d$ admits a finite volume uniformization by the unit ball $\bbB^2$ in $\bbC^2$. This produces one-cusped complex hyperbolic $2$-manifolds of arbitrarily large volume. As a consequence, the $3$-dimensional nilmanifold of Euler number $12d$ bounds geometrically for all odd $d \ge 1$.
\end{abstract}

\section{Introduction}\label{sec:Intro}

This introduction describes the main result of this paper along with some history and context in \S\ref{ssec:MainThm}, a more precise result in \S \ref{ssec:MainTech}, and an application to geometric bounding in \S\ref{ssec:Bound}.

\subsection{The main result, history, and context}\label{ssec:MainThm}

The purpose of this paper is to provide the first geometric construction of complete one-cusped complex hyperbolic $2$-manifolds of finite volume. The basic main result is as follows.

\begin{thm}\label{thm:Main}
For each odd $d \ge 1$, there is a one-cusped complex hyperbolic $2$-manifold of volume $16 \pi^2 d$.
\end{thm}

The existence of one-cusped manifolds was observed by the first author using Magma experimentation \cite[Thm.\ 1.4]{DerauxExp}, and a number of geometric properties were also recorded there, including the structure of the cusp cross-section. This paper gives an explicit, completely disjoint, computer-independent construction\footnote{While some computations are more easily done with the aid of computer algebra software, they can in principle all be done by hand.}. There is strong computational evidence suggesting that the examples in \cite{DerauxExp} are precisely those constructed here, and the inspiration for this paper was in fact similarities between the numerical invariants of those examples and the algebraic surfaces constructed in \cite{Polizzi}, but this paper does not address that question.

There is significant interest in existence of special metrics (e.g., K\"ahler--Einstein metrics) on a smooth projective variety $V$, or more generally metrics on $V \ssm D$ for a given divisor $D$; see for example \cite{Tian, TianYau}. Within this context, Theorem \ref{thm:Main} produces the first concrete example of a pair $(V, D)$ where $\dim(V) > 1$, $D$ is smooth and irreducible, and $V \ssm D$ admits a complete K\"ahler metric of constant biholomorphic sectional curvature $-1$ and finite volume. However, the primary interest in one-cusped complex hyperbolic manifolds stems from the fact that one-cusped locally symmetric manifolds of negative curvature are exceptional, rare, and elusive in (real) dimension greater than three.

With deep geometrization results in hand it becomes easy to construct $2$- and $3$-dimensional manifolds that admit complete, finite volume metrics of constant curvature $-1$ with exactly one cusp. Indeed, uniformization implies that any once-punctured Riemann surface of genus $g \ge 1$ admits such a metric, and Thurston's hyperbolization theorem \cite[Thm.\ 2.3]{ThurstonBulletin} leads to a plethora of examples (e.g., many knot complements in the $3$-sphere). However, shockingly little is known in higher dimensions, especially in the more general setting of locally symmetric manifolds with negative curvature (i.e., rank one). Hyperbolic $4$-manifolds with exactly one cusp were first constructed by Kolpakov and Martelli \cite[Thm.\ 1.1]{KolpakovMartelli}, and see \cite{Slavich, KolpakovSlavich1, KolpakovSlavich2, RatcliffeTschantz} for more examples. This paper and \cite[Thm.\ 1.4]{DerauxExp} provide a fourth known symmetric space of negative curvature admitting a one-cusped manifold quotient of finite volume; still no example of any kind is known in dimension greater than four.

For some perspective on why one-cusped manifolds are rare and special, the second author proved that the theory of arithmetic subgroups of algebraic groups is not as useful for producing one-cusped examples as one might think or hope. More specifically, for each $k \ge 1$ the rank one arithmetic locally symmetric spaces with $k$ cusps fall into finitely many commensurability classes over all possible universal covers and dimensions \cite[Thm.\ 1.1]{StoverCusp}. In particular, there is a dimension $d(n)$ so that $n$-cusped arithmetic locally symmetric manifolds of negative curvature cannot exist above dimension $d(n)$. For example, $1$-cusped arithmetic hyperbolic manifolds (in fact, orbifolds) cannot exist above dimension $30$ \cite[Thm.\ 1.3]{StoverCusp}. Thus it could well be the case that there are no one-cusped rank one locally symmetric spaces in sufficiently high dimensions.

\begin{rem}\label{rem:WithLuca}
The minimal possible volume of a complex hyperbolic $2$-manifold is $8 \pi^2 / 3$. Work of Kamishima \cite{Kamishima} asserts that if $M$ is a one-cusped finite volume complex hyperbolic $2$-manifold, then its cusp cross section must have trivial holonomy, which implies that $M$ admits a smooth toroidal compactification; see \S\ref{ssec:Bound} for further discussion. L. Di Cerbo and the second author proved that there is no minimal volume smooth toroidal compactification with one cusp \cite[Thm.\ 1.1]{DiCerboStover1}, so one-cusped manifolds must have volume at least $16 \pi^2 / 3$. On the other hand, Di Cerbo and the second author also proved that there is a two-cusped complex hyperbolic manifold realizing every possible volume \cite[Thm.\ 1.4]{DiCerboStover2}.
\end{rem}

\subsection{The main technical result}\label{ssec:MainTech}

Theorem \ref{thm:Main} is a direct consequence of the following more precise result.

\begin{thm}\label{thm:MainTech}
For each odd $d \ge 1$ there is a minimal smooth projective surface $Z_d$ of general type with $c_1^2(Z) = c_2(Z) = 6d$ and a smooth irreducible curve $E_d$ on $Z_d$ of genus one with self-intersection $-12d$ so that $Z_d \ssm E_d$ is uniformized by the unit ball $\bbB^2$ in $\bbC^2$.
\end{thm}

The pair $(Z_1, E_1)$ is constructed in \S\ref{sec:Construction} and the proof of Theorem \ref{thm:MainTech}, which immediately implies Theorem \ref{thm:Main} by Chern--Gauss--Bonnet, is the content of \S\ref{sec:Bootstrap}. The initial example $Z_1$ is one of the desingularized product-quotient surfaces with $p_g = q = 1$ studied by Polizzi \cite{Polizzi}. Briefly, there is a product $X = C_1 \times C_2$ of hyperbolic Riemann surfaces and a finite group $F$ acting on $X$ so that $Z_1$ is the minimal desingularization of $F \bs X$. The other examples are built using a covering construction.

\begin{rem}\label{rem:EqualHard}
As in the case of Riemann surfaces, the proof that $(Z_1, E_1)$ is a ball quotient applies a uniformization theorem, here due to Kobayashi \cite[Thm.\ 2]{Kobayashi} (also see \cite[Thm.\ 3.1]{TianYau}). A key difference makes dimension two and above much more subtle and special. A Riemann surface can be uniformized by a constant curvature metric if and only if a characteristic class (namely the Euler characteristic) has the appropriate sign. In higher dimensions, uniformization requires equality between Chern numbers, not merely an inequality, which is a significantly more stringent requirement.
\end{rem}

\subsection{Geometric bounding}\label{ssec:Bound}

One application of Theorem \ref{thm:MainTech} is to \emph{geometric bounding}. A famous theorem of Rohlin states that all closed, connected, orientable $3$-manifolds are diffeomorphic to the boundary of a compact $4$-manifold \cite{Rohlin}. For a $3$-manifold admitting one of the eight homogeneous model geometries \cite[\S 4]{Scott}, it is then of interest to know whether it can be realized as the boundary of a geometric $4$-manifold, and now there are obstructions, say from index theory. For example, Long and Reid studied when flat or hyperbolic manifolds can geometrically bound a hyperbolic manifold \cite{LongReid}. For flat manifolds, this means realizing the manifold as the cusp cross-section of a one-cusped hyperbolic manifold of one dimension higher, whereas in the hyperbolic case it means realizing the manifold as a totally geodesic boundary.

For infranil $3$-manifolds (i.e., those with Nil geometry), the relevant question is whether or not it can be realized as the cusp-cross section of a complex hyperbolic manifold. Restrictions on which manifolds can bound were given by unpublished work of Walter Neumann and Alan Reid, and work of Kamishima \cite{Kamishima} indicates that a manifold that bounds must have trivial holonomy. In other words, the cross-section should be a nilmanifold. In the language of toroidal compactifications, this is equivalent to saying that the complex hyperbolic manifold admits a \emph{smooth} toroidal compactification by an elliptic curve of self-intersection $-d$, where $d \ge 1$ and $d$ is the Euler number of the nilmanifold (as defined on \cite[p.\ 435]{Scott}). In the smooth toroidal case, Corollary \ref{cor:EulerRestrict} below shows that the Euler number is moreover divisible by four. The examples in this paper cover many, but not all, possibilities left after these reductions.

\begin{thm}\label{thm:Bound}
For every odd $d \ge 1$, the $3$-dimensional nilmanifold with Euler number $12 d$ geometrically bounds a complex hyperbolic manifold.
\end{thm}

As in the flat case \cite{LongReid2}, every infranil $3$-manifold is the cusp cross-section of some complex hyperbolic $2$-manifold, possibly with many cusps \cite{McReynolds, McReynolds2}. It would be interesting if every other nilmanifold with Euler number divisible by four can be realized as the cusp cross-section of a one-cusped manifold. It may even be the case that examples can be found that are commensurable with the examples in this paper.

\subsubsection*{Acknowledgments} This paper was written while Stover was a Professeur Invit\'e at Universit\'e Grenoble Alpes, and extends his profound thanks to the Institut Fourier for the excellent working environment. Stover was partially supported by Grants DMS-2203555 and DMS-2506896 from the National Science Foundation, along with award SFI-MPS-TSM-00014184 from the Simons Foundation. Revisions of this paper were completed while Stover was a CRM-Simons Professor at the Centre de Recherches Math\'ematiques, and he thanks both CRM and the Simons Foundation for their support. The authors thank the referee for a number of comments aimed at improving the paper.

\section{Background}\label{sec:Background}

This paper assumes basic familiarity with complex hyperbolic manifolds and the topology of smooth projective varieties. See for example \cite{Holzapfel, BPV} for basic references on what follows.

Let $\bbB^2$ be complex hyperbolic $2$-space, that is the unit ball in $\bbC^2$ with its metric of constant holomorphic sectional curvature $-1$. A complex hyperbolic manifold is a quotient $M = \Gam \bs \bbB^2$ of $\bbB^2$ by a torsion-free discrete group of holomorphic isometries. If $M$ has finite volume but is not compact, then it has a finite number of cusps, each diffeomorphic to the product of an infranil $3$-manifold $N$ with $[0, \infty)$, and $N$ is called a \emph{cusp cross-section} of $M$. See \cite[\S 2.2]{Dekimpe} for the definition of an infranil manifold.

When a cusp cross-section is a nilmanifold (i.e., has \emph{trivial holonomy}), the cusp is naturally diffeomorphic to a torus bundle over a punctured disk. Smoothly filling in the puncture gives a \emph{smooth toroidal compactification} of the cusp. Since it is the case relevant to this paper, suppose that $M$ has one cusp with smooth toroidal compactification $X$ obtained by adding an elliptic curve $E$. Then $(X, E)$ is called a \emph{ball quotient pair}. In that case, $X$ is in fact a smooth projective variety and $E$ is a smooth curve of genus one and self-intersection $-d$, where $d$ is the Euler number of the nil $3$-manifold in the cusp cross-section (see \cite[\S 4.2]{Holzapfel}).

Moreover, Kobayashi \cite[Thm.\ 2]{Kobayashi} proved a uniformization theorem that determines precisely when a pair $(X, E)$ consisting of a smooth projective variety and a smooth curve on $X$ determine a ball quotient pair. Suppose that $X$ is a smooth projective variety with canonical divisor $K_X$ and $E$ is a smooth curve of genus one on $X$. Then $X \ssm E$ admits a complete complex hyperbolic metric of finite volume if and only if the divisor $K_X + E$ is nef and big, there are no $(-2)$ curves on $X$ disjoint from $E$, and
\begin{equation}\label{eq:logBMY}
(K_X + E)^2 = 3 c_2(X),
\end{equation}
where $c_2(X)$ is the topological Euler characteristic. Here $(K_X + E)^2$ and $c_2(X)$ are the relative Chern numbers $c_1^2(X, E)$ and $c_2(X, E)$ of the pair $(X, E)$, respectively.

\begin{lem}\label{lem:IntersectionRestrict}
Let $(X, E)$ be a ball quotient pair with $E$ an elliptic curve of self-intersection $-d$. Then $d$ is divisible by four.
\end{lem}

\begin{pf}
Suppose $(X, E)$ is a ball quotient pair with $E$ an elliptic curve of self-intersection $-d$. The adjunction formula implies that
\begin{align*}
(K_X + E)^2 &= K_X^2 - E^2 \\
&= c_1^2(X) + d \\
&= 3 c_2(X)
\end{align*}
and so $d = 3 c_2(X) - c_1^2(X)$. However $c_1^2(X) = 12 \chi(\calO_X) - c_2(X)$ by Noether's formula, so
\[
d = 4 c_2(X) - 12 \chi(\calO_X),
\]
which is divisible by four.
\end{pf}

Reinterpreting Lemma \ref{lem:IntersectionRestrict} in terms of cusp cross-sections gives the following corollary.

\begin{cor}\label{cor:EulerRestrict}
If $M$ is a one-cusped complex hyperbolic manifold with cusp cross-section a nilmanifold $N$ with Euler number $d$, then $d$ is divisible by four.
\end{cor}

\section{Construction of the first example}\label{sec:Construction}

This section follows work of Polizzi \cite{Polizzi} to construct the crucial first example $(Z, E)$ needed to prove Theorem \ref{thm:Main}, and assumes that the reader is very familiar with the language of Fuchsian groups, for example as in \cite{Katok}. The notation $\Del(g; \frakn)$ will denote the Fuchsian group of signature $(g; \frakn)$. In other words, the quotient $\bbH^2 / \Del(g; \frakn)$ has genus $g$ and $a_j$ cone points of order $n_j$, where $\frakn = (n_1^{a_1}, \dots, n_k^{a_k})$ and any $n_j^1$ is simply given as $n_j$.

Briefly, the surface $Z$ is the minimal desingularization of the quotient of a product of curves under the action of a finite group. Curves $C_1$ and $C_2$ with action of the alternating group $\frakA_4$ are constructed in \S \ref{ssec:Genus4}. The quotient of $C_1 \times C_2$ by the diagonal action of $\frakA_4$ and its minimal smooth resolution $Z$ are studied in \S \ref{ssec:Product}. The elliptic curve $E$ on $Z$ is constructed in \S \ref{ssec:FindElliptic}. Finally, the fact that $(Z, E)$ is a ball quotient pair is proved in \S \ref{ssec:BallQ}.

\begin{rem}\label{rem:PolizziClass}
The example constructed in this section is of the kind covered by \cite[Prop.\ 7.3]{Polizzi}. However, the precise example presented here is not the one contained in Polizzi's proof. Polizzi is only concerned with showing a given triple consisting of two Fuchsian groups and a finite quotient produces at least one surface, and does not analyze all the possibilities for a given triple. Computer experiment indicates that each triple can in fact produce distinct surfaces that are not even homotopy equivalent.
\end{rem}

\subsection{The curves of genus four}\label{ssec:Genus4}

\subsubsection{The $(2,3,12)$-triangle group} \label{sssec:tg2312}

The construction begins with the $(2,3,12)$ triangle group
\[
\Gam_0 = \big\langle p,q,r\ |\ p^2,\, q^3,\, r^{12},\, pqr \big\rangle
\]
generated by three rotations by angles $\pi,2\pi/3,\pi/6$ in the vertices of a triangle
$T$ in the hyperbolic plane with respective angles $\pi/2,\pi/3,\pi/12$. Up to translation by an isometry there is a unique triangle in the hyperbolic plane with a given set of angles, which implies that there is a unique conjugacy class of discrete and faithful representations of $\Gam_0$ into $\PSL_2(\bbR)$ mapping $q$ and $r$ to anticlockwise rotations with fixed points at the appropriate vertices of $T$. We will frequently identify $\Gam_0$ with its image under such a representation.

Then
\[
F_0=T\cup\sigma(T)
\]
is a fundamental domain $F$ for the action of $\Gam_0$ on the hyperbolic plane, where $\sigma$ is a reflection in one of the sides of $T$. The images of $F_0$ under $\Gam_0$ tile the hyperbolic plane in such a way that there are $4$ (resp.\ $6$, $24$) copies of $T$ adjacent to the fixed point of $p$, (resp.\ $q$, $r$). Note also that every element of finite order in $\Gam_0$ is conjugate to a power of either $p$, $q$, or $r$.

\subsubsection{Two subgroups of index $6$ in $\Gam_0$}

Now consider the elements
\begin{align*}
a &= r^6 & g &= (prq)^{-1} \\
b &= q^2 & h &= (qpr)^{-1} \\
c &= (pr) q^2 (pr)^{-1} &\\
d &= (pr)^{-1} q^2 (pr) & &
\end{align*}
and let $\Gam_1$ (resp.\ $\Gam_2$) be the subgroup of $\Gam_0$ generated by $a,b,c,d$ (resp.\ $g,h$).

\begin{prop}\label{prop:g12}
Both $\Gam_1$ and $\Gam_2$ have index 6 in $\Gam_0$, and they have the following presentations in terms of the given generators:
  \begin{align*}
    \Gam_1 &= \big\langle a,b,c,d \big\rangle \\
           &\cong \big\langle a,b,c,d\ |\ a^2,\, b^3,\, c^3,\, d^3,\, abcd \big\rangle \\
    \Gam_2 &= \big\langle g,h \big\rangle \\
           &\cong \big\langle g,h\ |\ [g,h]^2 \big\rangle
  \end{align*}
\end{prop}

\begin{figure}[htbp]
\centering
\includegraphics[width=0.48\textwidth]{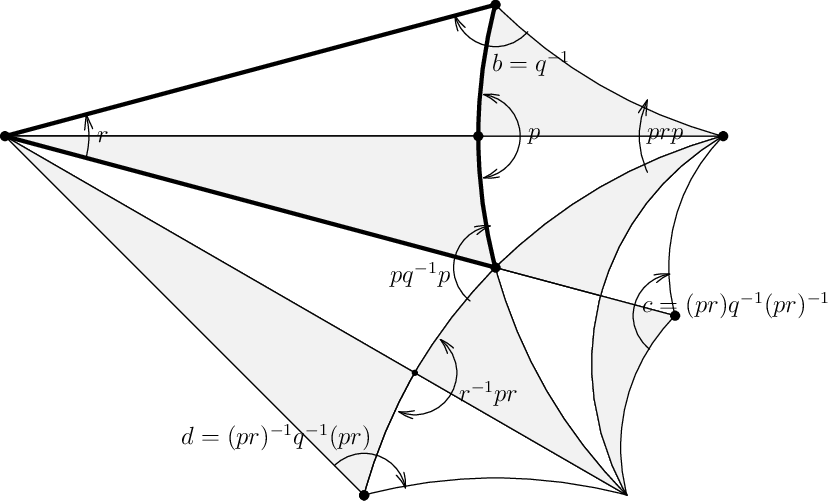}\quad
\includegraphics[width=0.48\textwidth]{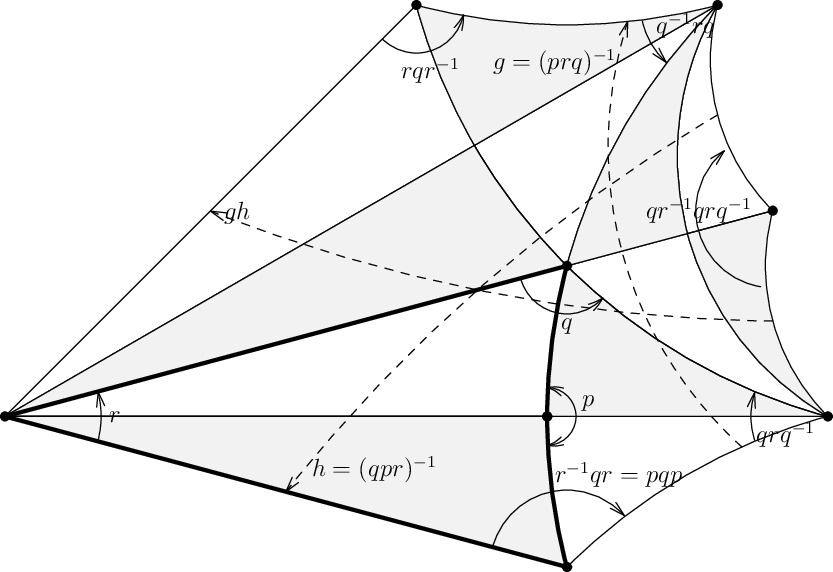}
\caption{Fundamental domains for $\Gamma_1$ (left) and $\Gamma_2$ (right), as copies of a fundamental domain for $\Gamma_0$ (shown in bold on both pictures)}\label{fig:fundom}
\end{figure}

\begin{pf}
The fact that $a^2=b^3=c^3=d^3=abcd=[g,h]^2=\Id$ in $\Gam_0$ follows from direct computation. For example, note that $q$ and $pr=r^{-1}q^{-1}r$ both have order $3$. Moreover, $p=qr=r^{-1}q^{-1}$ since $p = p^{-1}$, so
\begin{align*}
abcd&=r^6\cdot q^2\cdot(pr) q^2(pr)^{-1}\cdot (pr)^{-1}q^2(pr)\\
&= r^6 q^{-1} q r^2 q^2 r^{-1} q^{-1} r q^2 q r^2\\
&= r^{12}\\
&= \Id
\end{align*}
A similar computation shows that the commutator $[g,h]=ghg^{-1}h^{-1}$ has order 2.

The Poincar{\'e} polygon theorem is now used to show that these relations give presentations for $\Gam_1$ and $\Gam_2$. The pictures in Figure~\ref{fig:fundom} describe explicit polygons obtained by taking images of a standard fundamental domain $F$ for $\Gam_0$, with side pairing identifications given by $a,b,c,d$ to the left and $g,h,gh$ to the right. Explicitly,
\begin{align*}
F_1 &= F_0 \cup r^{-1}(F_0) \cup p(F_0) \cup r^{-1}pr(F_0) \cup pqpr^{-1}(F_0) \cup pr(F_0) \\
F_2 &= F_0 \cup p(F_0) \cup q(F_0) \cup r(F_0) \cup rp(F_0) \cup qp(F_0)
\end{align*}
are fundamental domains for $\Gam_1$ and $\Gam_2$.

We show one sample computation showing that the relevant maps are indeed side-pairing identifications, namely that ${g=(prq)^{-1}}$ maps the geodesic segment joining the fixed points of $x=pqp=r^{-1}qr$ and $y=qrq^{-1}$ to the geodesic segment joining the fixed points of $x^\prime=rqr^{-1}$ and $y^\prime=q^{-1}rq$. This follows from the fact that $gxg^{-1}=x^\prime$ and $gyg^{-1}=y^\prime$. To justify $gxg^{-1}=x^\prime$,
\begin{align*}
gxg^{-1}&=q^{-1}r^{-1}p\cdot r^{-1}qr\cdot prq\\
&= q^{-1}r^{-1}\cdot q\cdot rq = q^{-1}r^{-1}q^{-1}qqrq = rqr^{-1}=x^\prime
\end{align*}
as claimed. The other calculations are of the same level of difficulty.
  
To verify the hypotheses of the Poincar{\'e} polygon theorem we must determine cycle transformations associated with each vertex of the polygon and check that each gives a rotation by an angle of the form $2\pi/k$ for some $k\in\bbN$. One check is done, and the rest are left to the reader. For the left part of Figure~\ref{fig:fundom}, the only cycle transformation that is not evident from the image is the one associated with the fixed point of $r$. Tracking the identifications given by $a,b,c$ and $d$, the cycle transformation is
\begin{align*}
bcd &= q^{-1}\cdot pr q^{-1}(pr)^{-1}\cdot(pr)^{-1}q^{-1}pr\\
&= q^{-1}qrrq^{-1}r^{-1}qr r^{-1}qrq^{-1}qrr\\
&= r^2q^{-1}r^{-1}q^{-1}r^3\\
&= r^2q^{-1}qr r^3 \\
&= r^6
\end{align*}
which is a rotation by $\pi$.
\end{pf}

%
%
%

\begin{figure}[h]
\centering
\begin{tikzpicture}[scale=0.5]
\clip (8,8) rectangle (-8,-6.75);
\draw (0,5) ellipse (2 and 2);
\begin{scope}
\clip (2,5) rectangle (-2, 3);
\draw (0,5) ellipse (2 and 0.5);
\end{scope}
\begin{scope}
\clip (2,5) rectangle (-2, 7);
\draw [dashed] (0,5) ellipse (2 and 0.5);
\end{scope}
\draw (-5,0) ellipse (2 and 2);
\begin{scope}
\clip (-7,0) rectangle (-3, -2);
\draw (-5,0) ellipse (2 and 0.5);
\end{scope}
\begin{scope}
\clip (-7,0) rectangle (-3, 2);
\draw[dashed] (-5,0) ellipse (2 and 0.5);
\end{scope}
\draw (5,0) ellipse (3 and 1.5);
\begin{scope}
\clip (5,-1.8) ellipse (3 and 2.5);
\draw (5,2.2) ellipse (3 and 2.5);
\end{scope}
\begin{scope}
\clip (5,2.2) ellipse (3 and 2.5);
\draw (5,-2.2) ellipse (3 and 2.5);
\end{scope}
\draw (0,-5) ellipse (3 and 1.5);
\begin{scope}
\clip (0,-6.8) ellipse (3 and 2.5);
\draw (0,-2.8) ellipse (3 and 2.5);
\end{scope}
\begin{scope}
\clip (0,-2.8) ellipse (3 and 2.5);
\draw (0,-7.2) ellipse (3 and 2.5);
\end{scope}
\draw[->, shorten <=1.25cm, shorten >=1.25cm] (5,0) -- node [above right] {\footnotesize{$\bbZ/6$}} (0,5);
\draw[->, shorten <=1.25cm, shorten >=1.25cm] (-5,0) -- node [above left] {\footnotesize{$6$-to-$1$}} (0,5);
\draw[->, shorten <=1.25cm, shorten >=1.25cm] (0,-5) -- node [below right] {\footnotesize{$\bbZ/3$}} (5,0);
\draw[->, shorten <=1.25cm, shorten >=1.25cm] (0,-5) -- node [below left] {\footnotesize{$\bbZ/3$}} (-5,0);
\filldraw[red] (0,7) circle (3pt) node [above] {\footnotesize{$12$}};
\filldraw (-2,5) circle (3pt) node [left] {\footnotesize{$2$}};
\filldraw[blue] (2,5) circle (3pt) node [right] {\footnotesize{$3$}};
\filldraw[red] (-5,2) circle (3pt) node [above] {\footnotesize{$2$}};
\filldraw[blue] (-7,0) circle (3pt) node [left] {\footnotesize{$3$}};
\filldraw[blue] (-3,0) circle (3pt) node [right] {\footnotesize{$3$}};
\filldraw[blue] (-5,-0.5) circle (3pt) node [below] {\footnotesize{$3$}};
\filldraw[red] (7.25,0) circle (3pt) node [above] {\footnotesize{$2$}};
\filldraw[red] (2.25,-5) circle (3pt) node [above] {\footnotesize{$2$}};
\filldraw[red] (-2.25,-5) circle (3pt) node [above] {\footnotesize{$2$}};
\filldraw[red] (0,-6.25) circle (3pt) node [above] {\footnotesize{$2$}};
\end{tikzpicture}
\caption{Coverings of hyperbolic orbifolds}\label{fig:Covers1}
\end{figure}
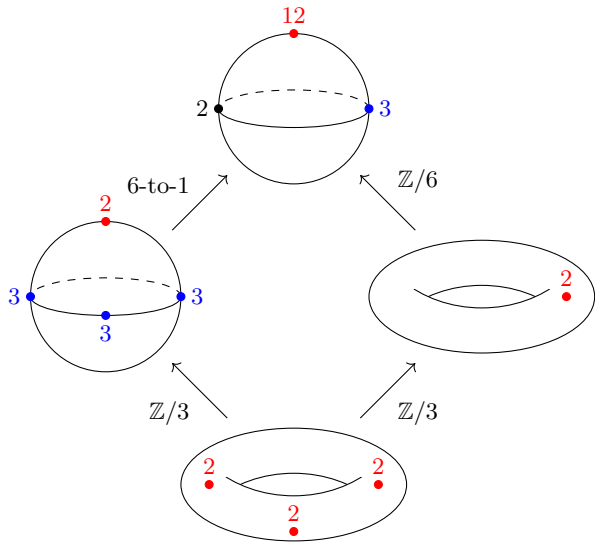

In terms of orbifold quotients, there are maps of orbifolds depicted in the top half of Figure \ref{fig:Covers1}. In particular, $\Gam_1 \cong \Del(0; 2, 3^3)$ and $\Gam_2 \cong \Del(1; 2)$.

\medskip

Now consider the homomorphisms $\rho_j : \Gam_j \to \frakA_4$ induced
by
\begin{align*}
\rho_1(a) &= (1\, 3)(2\, 4) & \rho_2(g) &= (1\, 2\, 3) \\
\rho_1(b) &= (1\, 3\, 2) & \rho_2(h) &= (1\, 4\, 2) \\
\rho_1(c) &= (1\, 2\, 4) & & \\
\rho_1(d) &= (1\, 2\, 4) & &
\end{align*}
with the convention that the alternating group $\frakA_4$ acts on $\{1,2,3,4\}$ on the right\footnote{While not the authors' preference, since actions in this paper are on the left, this choice is for consistency with standard computer algebra software.}. Let $\Lam_j$ be the kernel of $\rho_j$. It is easy to check that $\rho_j$ is surjective, so $\Lam_j$ has index 12 in $\Gam_j$. Moreover every element of finite order in $\Gam_1$ (resp.\ $\Gam_2$) is conjugate to a power of $a$, $b$, $c$ or $d$ (resp.\ a power of $[g,h]$). Using this, one easily checks for $j=1$ and $2$ that the homomorphism $\rho_j$ restricts to an injection on every finite subgroup of $\Gam_j$, so $\Lam_j$ is torsion-free. Using a fixed identification of $\Gam_0$ with a subgroup of $\PSL_2(\bbR)$ as in Section~\ref{sssec:tg2312} and its restriction to the subgroups $C_j$, we obtain a Riemann surface $C_j = \Lam_j \bs \bbH^2$ with an action of $\frakA_4$ induced by $\rho_j$.

The orbifold Euler characteristic of $\Gam_0 \bs \bbH^2$ is
\[
\chi^{orb}(\Gam_0 \bs \bbH^2)=-1+\frac{1}{2}+\frac{1}{3}+\frac{1}{12}=-\frac{1}{12},
\]
so $\chi^{orb}(\Gam_j \bs \bbH^2)=6\chi^{orb}(\Gam_0 \bs \bbH^2)=-1/2$ and $\chi(C_j)=-6$. It follows that $C_1$ and $C_2$ both have genus four. The following lemma gives more information concerning the $\frakA_4$ action on each $C_j$.

\begin{lem}\label{lem:CjAut}
Each order two element of $\frakA_4$ has exactly two fixed points on $C_1$. These fixed points are precisely the six lifts to $C_1$ of the unique cone point of order two on $\Gam_1 \bs \bbH^2$, hence the $\frakA_4$ action on these six points is transitive. No order three element of $\frakA_4$ has fixed point on $C_2$.
\end{lem}

\begin{pf}
The proof is standard branched covering theory. The case of $C_1$ is explained and $C_2$ is left to the reader. First, fixed points on $C_1$ must map to cone points on the quotient orbifold $\Gam_1 \bs \bbH^2$. The unique conjugacy class of order $2$ elements of $\frakA_4$, which consists of the images of order two elements of $\Gam_1$, is then associated with the unique cone point of order two on $\frakA_4 \bs C_1$, which has $12 / 2 = 6$ lifts to $C_1$ with transitive $\frakA_4$ action. Since conjugate elements must have the same number of fixed points, each must fix exactly two of the six preimages.
\end{pf}

\subsection{The product action}\label{ssec:Product}

Retaining the notation of \S\ref{ssec:Genus4}, define $X = C_1 \times C_2$ and consider the diagonal action of $\frakA_4$ on $X$. Lemma \ref{lem:CjAut} implies that each $3$-cycle acts freely on $X$, since it acts freely on the second coordinate. Each order two element has exactly four fixed points, namely the points of the form $(z, w)$ with $z$ one of its two fixed points on $C_1$ and $w$ one of its two fixed points on $C_2$.

This gives a total of twelve points on the product on which the action of $\frakA_4$ is faithful by Lemma \ref{lem:CjAut}, hence
\[
Y = \frakA_4 \bs X
\]
has exactly two singular points, which are quotient singularities of type $\mathrm{A}_1$. Indeed, the action of an order two element $\sig \in \frakA_4$ on the tangent space to a fixed point $(z,w)$ is by $-\Id$, since the action on each individual coordinate is locally a rotation of order two. Let $\varphi : Z \to Y$ be the
minimal resolution of singularities of $Y$, which gives a diagram
\[
\begin{tikzcd} X \arrow[d, "\pi"] & \\ Y & \arrow[l, "\varphi"] Z \end{tikzcd}
\]
where $\pi$ is projection for the $\frakA_4$ action. The minimal resolution has the property that if $y_1, y_2$ are the singular points of $Y$, then $\varphi$ is an isomorphism on $Z \ssm \varphi^{-1}(\{y_1, y_2\})$ and each $\varphi^{-1}(y_j)$ is a smooth rational curve $F_j$ of self-intersection $-2$. See \cite[\S III.1-7]{BPV}.

Note that $X$ has Euler characteristic $36$, the action of $\frakA_4$ has twelve fixed points, and each singular point of $Y$ is replaced in $Z$ by $\bbP^1$ (i.e., topologically a $2$-sphere), hence the Euler characteristic of $Z$ is
\[
c_2(Z) = \frac{1}{12}(36 - 12) + 2 \times 2 = 6.
\]
Moreover, $K_Z = \varphi^* K_Y$ since the resolution is of two $\mathrm{A}_1$ singularities \cite[\S I.1]{Reid}. Indeed, $K_Z$ can be written as $\varphi^* K_Y + \al_1 F_1 + \al_2 F_2$ for some $\al_1, \al_2\in\bbQ$
by basic properties of the canonical bundle, but adjunction implies
that each $\al_j$ must be zero. Therefore
\[
c_1(Z)^2 = \varphi^*(K_Y)^2 = \frac{1}{12} K_X^2 = 6,
\]
since $c_1^2(X) = 2 c_2(X)$ for smooth compact quotients of $\bbH^2 \times \bbH^2$ by Hirzebruch proportionality \cite[Satz 2]{Hirzebruch}. Thus $Z$ is a smooth surface such that ${c_1^2(Z) = c_2(Z) = 6}$.

Finally, observe that $Z$ is minimal of general type; see \cite[Prop.\ 5.6]{Polizzi}. Moreover, as discussed in \cite[\S 5]{Polizzi}, $p_g(Z) = q(Z) = 1$ where the Albanese map of $Z$ is $\varphi$ followed by the projection of $Y$ onto $\Gam_2 \bs \bbH^2$ induced by projection of $C_1 \times C_2$ onto the second factor. The contents of this subsection are collected in the following result.

\begin{prop}\label{prop:DescribeQuo}
Let $C_1$ and $C_2$ be the genus four curves with $\frakA_4$ action described in \S \ref{ssec:Genus4}. Set $X = C_1 \times C_2$ and consider the diagonal action of $\frakA_4$ with quotient $Y$. Then $Y$ has two singularities of type $\mathrm{A}_1$ and its minimal resolution $Z$ is a minimal surface of general type with $p_g = q = 1$ and $c_1^2 = c_2 = 6$. The Albanese map of $Z$ is induced by the natural projection of $Y$ onto the smooth curve $\frakA_4 \bs C_2$ of genus one.
\end{prop}

\subsection{The elliptic curve $E$}\label{ssec:FindElliptic}

The surface $X$ in \S \ref{ssec:Product} is uniformized by $\Lam_1 \times \Lam_2$, where the action on $\bbH^2 \times \bbH^2$ is induced by the product action of $\Gam_0 \times \Gam_0$. Consider the intermediate group $\Gam_1 \times \Gam_2$, which produces a sequence of orbifold covers
\[
\begin{tikzcd}
\Gam_1 \bs \bbH^2 \times \Gam_2 \bs \bbH^2 \\
Y \arrow[u] \\
C_1 \times C_2 \arrow[u, "\frakA_4"] \arrow[uu, bend right=45, "\frakA_4 \times \frakA_4" right]
\end{tikzcd}
\]
where the bottom covering is the diagonal action of $\frakA_4$ and the composition of the two covers is the product action of $\frakA_4 \times \frakA_4$. The goal of this subsection is to use this diagram to prove that the image of the diagonal $\wt{D}$ of $\bbH^2 \times \bbH^2$ in $Y$ is a singular curve whose normalization has genus one.

The stabilizer in $\Gam_1 \times \Gam_2$ of $\wt{D}$ is naturally isomorphic to the subgroup $\Gam_1 \cap \Gam_2$ of $\Gam_0$. The relations
\begin{align*}
  t_1 &= g^3 & s_1 &= [g,h] \\
&= b d c & &= a \\
t_2 &= g^{-1} h & s_2 &= (h g h^{-1}) [g,h] (h g h^{-1})^{-1} \\
 &= c^{-1} b & &= (d^2 c b^2) a (d^2 c b^2)^{-1} \\
& & s_3 &= (h g^{-1} h^{-1}) [g, h] (h g^{-1} h^{-1})^{-1} \\
& & &= d a d^{-1}
\end{align*}
in $\Gam_0$ imply that $\Gam_1 \cap \Gam_2$ contains the group $\Gam_3$ generated by these elements. A long (but not particularly difficult) calculation allows one to verify the relation
\[
[t_1, t_2] = s_1 s_2 s_3,
\]
so $\Gam_3$ is isomorphic to a quotient of
\[
\Gam(1; 2^3)\cong \left\langle s_1,s_2,s_3,t_1,t_2\ |\ s_1^2=s_2^2=s_3^2=\Id, [t_1,t_2]=s_1s_2s_3\right\rangle.
\]
We will soon see that it is actually isomorphic to $\Gam(1;2^3)$ by computing its index in $\Gam_j$.

Note that there are homomorphisms $\psi_j:\Gam_j\rightarrow\bbZ / 3=\langle \sig \rangle$, uniquely characterized by the fact that $\psi_1$ sends $a$ to the identity and $b,c,d$ to $\sig$ and $\psi_2$ sends both $g$ and $h$ to $\sig$. Set $K_j=\mathrm{ker}(\psi_j)$. From the above relations, it is clear that $\Gam_3\subset K_1$, and $\Gam_3\subset K_2$. Moreover,
\[
\chi^{orb}(K_1 \bs \bbH^2)=\chi^{orb}(K_2 \bs \bbH^2)=3\chi^{orb}(\Gam_j \bs \bbH^2)=-\frac{3}{2}=\chi^{orb}(\Gam(1;2^3) \bs \bbH^2),
\]
which shows that $K_1=K_2=\Gam_3$. It also follows that $\Gam_3$ is normal in both $\Gam_1$ and $\Gam_2$ and that $\Gam_3 \bs \bbH^2$ is therefore the orbifold of genus one with three cone points of order two completing the diagram in Figure \ref{fig:Covers1}.

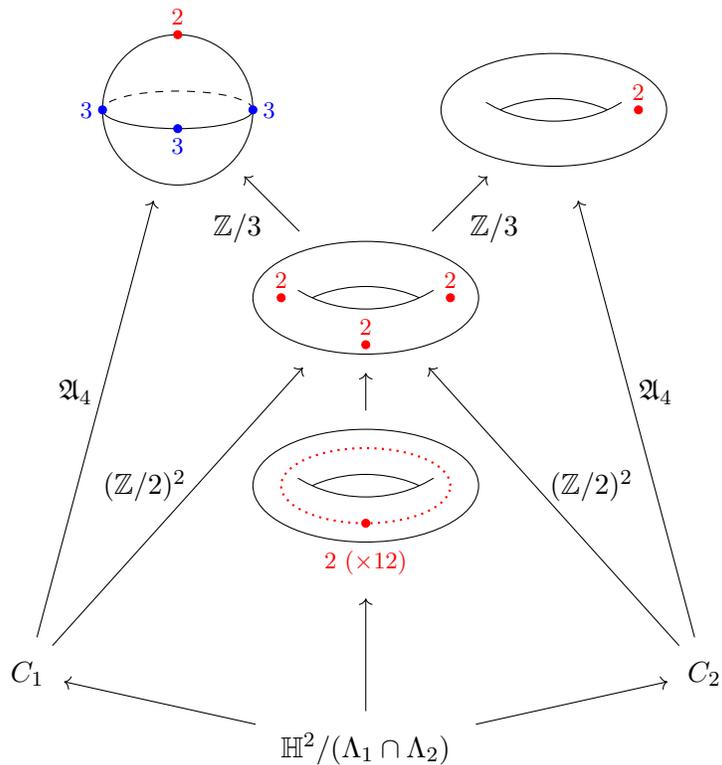
\begin{figure}[h]
\centering
\begin{tikzpicture}[scale=0.5]
\draw (-5,0) ellipse (2 and 2);
\begin{scope}
\clip (-7,0) rectangle (-3, -2);
\draw (-5,0) ellipse (2 and 0.5);
\end{scope}
\begin{scope}
\clip (-7,0) rectangle (-3, 2);
\draw[dashed] (-5,0) ellipse (2 and 0.5);
\end{scope}
\draw (5,0) ellipse (3 and 1.5);
\begin{scope}
\clip (5,-1.8) ellipse (3 and 2.5);
\draw (5,2.2) ellipse (3 and 2.5);
\end{scope}
\begin{scope}
\clip (5,2.2) ellipse (3 and 2.5);
\draw (5,-2.2) ellipse (3 and 2.5);
\end{scope}
\draw (0,-5) ellipse (3 and 1.5);
\begin{scope}
\clip (0,-6.8) ellipse (3 and 2.5);
\draw (0,-2.8) ellipse (3 and 2.5);
\end{scope}
\begin{scope}
\clip (0,-2.8) ellipse (3 and 2.5);
\draw (0,-7.2) ellipse (3 and 2.5);
\end{scope}
\draw (0,-10) ellipse (3 and 1.5);
\begin{scope}
\clip (0,-11.8) ellipse (3 and 2.5);
\draw (0,-7.8) ellipse (3 and 2.5);
\end{scope}
\begin{scope}
\clip (0,-7.8) ellipse (3 and 2.5);
\draw (0,-12.2) ellipse (3 and 2.5);
\end{scope}
\node at (-9, -15) {$C_1$};
\node at (9, -15) {$C_2$};
\node at (0, -17) {$\bbH^2 / (\Lam_1 \cap \Lam_2)$};
\draw[->, shorten <=1.25cm, shorten >=1.25cm] (0,-5) -- node [below right] {$\bbZ/3$} (5,0);
\draw[->, shorten <=1.25cm, shorten >=1.25cm] (0,-5) -- node [below left] {$\bbZ/3$} (-5,0);
\draw[->, shorten <=1cm, shorten >=1cm] (0,-10) -- (0,-5);
\draw[->, shorten <=0.5cm, shorten >=1.25cm] (-9,-15) -- node [left] {$\frakA_4$} (-5,0);
\draw[->, shorten <=0.5cm, shorten >=1.25cm] (-9,-15) -- node [left] {$(\bbZ / 2)^2$} (0,-5);
\draw[->, shorten <=0.5cm, shorten >=1.25cm] (9,-15) -- node [right] {$\frakA_4$} (5,0);
\draw[->, shorten <=0.5cm, shorten >=1.25cm] (9,-15) -- node [xshift=0.75cm] {$(\bbZ / 2)^2$} (0,-5);
\draw[->, shorten <=1.5cm, shorten >=0.5cm] (0,-17) -- (-9,-15);
\draw[->, shorten <=1.5cm, shorten >=0.5cm] (0,-17) -- (9,-15);
\draw[->, shorten <=0.5cm, shorten >=0.5cm] (0,-17) -- (0,-12);
\filldraw[red] (-5,2) circle (3pt) node [above] {\footnotesize{$2$}};
\filldraw[blue] (-7,0) circle (3pt) node [left] {\footnotesize{$3$}};
\filldraw[blue] (-3,0) circle (3pt) node [right] {\footnotesize{$3$}};
\filldraw[blue] (-5,-0.5) circle (3pt) node [below] {\footnotesize{$3$}};
\filldraw[red] (7.25,0) circle (3pt) node [above] {\footnotesize{$2$}};
\filldraw[red] (2.25,-5) circle (3pt) node [above] {\footnotesize{$2$}};
\filldraw[red] (-2.25,-5) circle (3pt) node [above] {\footnotesize{$2$}};
\filldraw[red] (0,-6.25) circle (3pt) node [above] {\footnotesize{$2$}};
\filldraw[red] (0,-11) circle (3pt) node [yshift=-15pt] {\footnotesize{$2\ (\times 12)$}};
\draw[red, dotted, thick] (0,-10) ellipse (2.25 and 1);
\end{tikzpicture}
\caption{Further coverings of hyperbolic orbifolds}\label{fig:Covers2}
\end{figure}

If $\Om < \Gam_1 \times \Gam_2$ is the subgroup associated with $Y$, then $\Om$ is the preimage in $\Gam_1 \times \Gam_2$ of the diagonal subgroup of $\frakA_4 \times \frakA_4$ under $\rho_1 \times \rho_2$. Thus the stabilizer in $\Om$ of $\wt{D}$ is
\[
\Gam_4 = \big\{\gam \in \Gam_3\ :\ \rho_1(\gam) = \rho_2(\gam) \big\}.
\]
Direct checks show that
\begin{align*}
\rho_1(t_1) &= (1\, 3)(2\, 4) & \rho_2(t_1) &= \Id \\
\rho_1(t_2) &= (1\, 4)(2\, 3) & \rho_2(t_2) &= (1\, 3)(2\, 4)
\end{align*}
and that $\rho_1$ and $\rho_2$ agree on each $s_j$. Then
\begin{align*}
\rho_1(t_1^2) &= \rho_2(t_1^2) \\
\rho_1(t_2^2) &= \rho_2(t_2^2)
\end{align*}
and so $\Gam_4$ contains the index four subgroup of $\Gam_3$ isomorphic to $\Del(1; 2^{12})$ induced by the unique $(\bbZ / 2)^2$ unramified cover of the torus. From the fact that $\rho_1$ and $\rho_2$ differ on $t_j^{\ell_1} t_k^{\ell_2}$ with $\ell_1, \ell_2 \in \{0,1\}$ not both zero, it follows that $\Gam_4 \cong \Del(1; 2^{12})$ with natural generating set $t_1^2$, $t_2^2$, and all the appropriate conjugates of each $s_j$. For instance, one can take $t_1s_jt_1^{-1}$, $t_2s_jt_2^{-1}$, $t_1t_2s_j(t_1t_2)^{-1}$ for $j=1,2,3$. This leads to the diagram of orbifold coverings depicted in Figure \ref{fig:Covers2}.

By construction, the groups defined in this section give a commutative diagram of (horizontal) immersions
\begin{equation}\label{eq:Immersions}
\begin{tikzcd}
\Gam_0 \bs \bbH^2 \arrow[r, hookrightarrow] & \Gam_0 \bs \bbH^2 \times \Gam_0 \bs \bbH^2 \\
\Gam_3 \bs \bbH^2 \arrow[r, loop-math to] \arrow[u] & \Gam_1 \bs \bbH^2 \times \Gam_2 \bs \bbH^2 \arrow[u] \\
\Gam_4 \bs \bbH^2 \arrow[r, loop-math to] \arrow[u] & Y \arrow[u]
\end{tikzcd}
\end{equation}
associated with projections of the diagonal $\wt{D}$ of $\bbH^2 \times \bbH^2$. While the diagonal of $\Gam_0 \bs \bbH^2 \times \Gam_0 \bs \bbH^2$ is certainly embedded, the coverings under consideration are orbifold covers, hence it does not follow that $\Gam_j \bs \bbH^2$ embeds for $j = 3,4$. In fact, this is not the case.

\begin{prop}\label{prop:SelfIntersection}
The image $\wh{E}$ of $\Gam_4 \bs \bbH^2$ in $Y$ meets each of the two $\mathrm{A}_1$ singularities with multiplicity six. Away from those points, the map is an embedding. Consequently, the proper transform $E$ of $\wh{E}$ to the resolution $Z$ of $Y$ is a smooth curve of genus one with self-intersection $-12$.
\end{prop}

The proof of Proposition \ref{prop:SelfIntersection} is quite involved, and requires some preliminary lemmas (Lemmas~\ref{lem:SIstep1} through~\ref{lem:SIstep3}). Notation established in each step will be used freely in each subsequent step. For $n \in\{2, 3,12\}$, let $\wt{z}_n \in \wt{D}$ be the point fixed by the diagonal action of $p$, $q$, and $r$ on $\bbH^2 \times \bbH^2$, respectively. Then the image $z_n$ of $\wt{z}_n$ on $\Gam_0 \bs \wt{D}$ is its cone point of order $n$.

\begin{lem}\label{lem:SIstep1}
All intersections of $\wt{D}$ with its orbit under $\Gam_0 \times \Gam_0$ arise from translates under the diagonal action of $\Gam_0$ of the three configurations depicted in Figure \ref{fig:Diag0}.
\end{lem}

\begin{pf}
Consider the diagonal embedding of $D_0 = \bbH^2 / \Gam_0$ in
\[
X_0 = \Gam_0 \bs \bbH^2 \times \Gam_0 \bs \bbH^2 \simeq \bbP^1 \times \bbP^1.
\]
At a cone point $z_n$ on $D_0$ of order $n \in \{2,3,12\}$ with lift $\wt{z}_n \in \wt{D}$ as above, the orbifold covering by $\bbH^2 \times \bbH^2$ has local group $(\bbZ / n)^2$ and there are $n$ preimages of $\Gam_0 \bs \bbH^2$ through $\wt{z}_n$, namely the graphs of the various powers of the associated elliptic element of $\PSL_2(\bbR)$. This gives the configurations of isometric embeddings of $\bbH^2$ in the product depicted in Figure \ref{fig:Diag0}.

\begin{figure}[h]
\centering
\begin{tikzpicture}
\draw (-3,-3) -- (3,3) node [above right] {$\wt{D}$};
\draw (-2.5,-1.5) -- (-1.5,-2.5) node [xshift=0.25cm, yshift=-0.25cm] {\footnotesize{$\big\{(z, p(z))\big\}$}};
\filldraw[black] (-2,-2) circle (1.25pt) node [xshift=-0.5cm] {$\wt{z}_2$};
\draw (0,0.5) -- (0,-0.5) node [yshift=-0.25cm, xshift=0.25cm] {\footnotesize{$\big\{(z, q^2(z))\big\}$}};
\draw (0.5,0) -- (-0.5,0) node [xshift=1.8cm] {\footnotesize{$\big\{(z, q(z))\big\}$}};
\filldraw[blue] (0,0) circle (1.25pt) node [above left] {$\wt{z}_3$};
\draw (2,2.5) -- (2,1.5);
\draw (2.5,2) -- (1.5,2) node [xshift=2cm] {\footnotesize{$\big\{(z, r^j(z))\big\}$}};
\filldraw[red] (2,2) circle (1.25pt) node [above left] {$\wt{z}_{12}$};
\draw[thick, densely dotted] (2.25, 2) arc (0:-90:0.25);
\end{tikzpicture}
\caption{The $\Gam_0 \times \Gam_0$ orbits of $\wt{D}$ that meet $\wt{D}$}\label{fig:Diag0}
\end{figure}
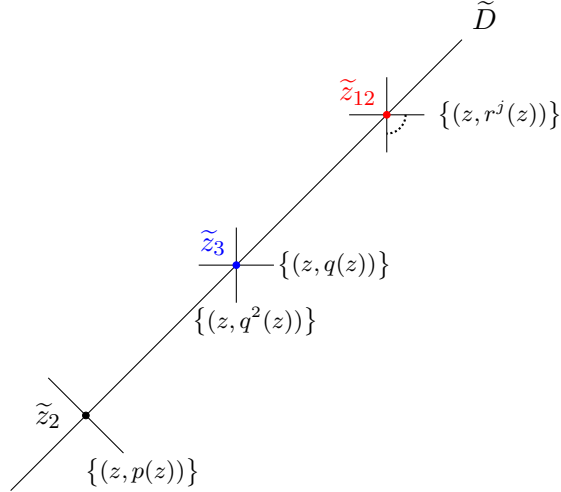

Moreover, if $(\al, \beta) \in \Gam_0 \times \Gam_0$ with $\al \neq \beta$ and $(\al, \beta)(\wt{D}) \cap \wt{D}$ is nonempty, then the intersection is a point
\[
(w,w) = (\al(z), \beta(z)),
\]
and thus $\beta^{-1} \al(z) = z$ and so $\beta^{-1} \al$ is a torsion element of $\Gam_0$. Since the conjugacy classes of torsion elements of $\Gam_0$ are represented by powers of conjugates of $p,q,r$, it follows that $(z,z)$ is a $\Gam_0$-translate of some $\wt{z}_j$. Thus all intersections of $\wt{D}$ with its orbit under $\Gam_0 \times \Gam_0$ arise from translates of the three configurations depicted in Figure \ref{fig:Diag0} under the diagonal action of $\Gam_0$. This proves the lemma.
\end{pf}

Next, consider $X_3 = \Gam_1 \bs \bbH^2 \times \Gam_2 \bs \bbH^2$ and $D_3 = \Gam_3 \bs \bbH^2$ from the middle row of Equation \eqref{eq:Immersions}, recalling that $D_3$ is genus one with three cone points of order two.

\begin{lem}\label{lem:SIstep2}
The immersion $D_3 \looparrowright X_3$ is an embedding away from the three cone points of $D_3$ and the cone points all map to a single point of $X_3$, as depicted in Figure \ref{fig:EmbedD3}.
\end{lem}

\begin{pf}
The map from $X_3$ to $X_0$ has degree $36$, since $[\Gam_0 : \Gam_j] = 6$ for $j = 1,2$. The map from $D_3$ to $D_0$ has degree $18$; see Figure~\ref{fig:Covers2}. It follows that the preimage of $D_0 \subset X_0$ in $X_3$ has two irreducible components, namely the image of $D_3$ and another curve $D_3^\prime$, where the map $D_3 \looparrowright X_3$ is induced by the orbifold covering projection $\Gam_3 \bs \bbH^2 \to \Gam_j \bs \bbH^2$ in each coordinate $j = 1,2$. Since the $\Gam_0 \times \Gam_0$ orbit of $\wt{D}$ only intersects $\wt{D}$ nontrivially at lifts of orbifold points of $D_0$ by Lemma \ref{lem:SIstep1}, the irreducible components of the preimage of $D_0$ and their intersection properties can be understood by understanding the preimages in $X_3$ of the cone points of $D_0$.
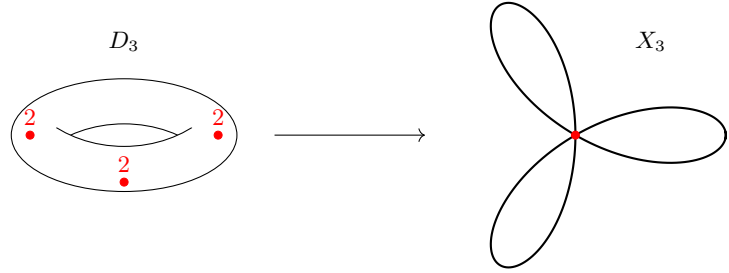
\begin{figure}[h]
\centering
\begin{tikzpicture}[scale=0.5]
\draw[variable=\t, domain=-1.6:1.6, samples=250]
plot ({(1/2)*(1+cos(\t 3)^2)*(4*(4*sin(\t r)^4-3*sin(\t r)^2))}, {(1/2)*(1+cos(\t 3)^2)*(4*(-sin(\t r)*cos(\t r)*(4*sin(\t r)^2-3)))});
\draw[variable=\t, domain=-1.6:1.6, samples=250]
plot ({(1/3)*(1+cos(\t 3)^2)*(4*(4*sin(\t r)^4-3*sin(\t r)^2))}, {(1/3)*(1+cos(\t 3)^2)*(4*(-sin(\t r)*cos(\t r)*(4*sin(\t r)^2-3)))});
\draw (-2,3.4641) to [bend left=40] (-1.33333,2.3094);
\draw[dashed] (-2,3.4641) to [bend right=40] (-1.33333,2.3094);
\draw (-1.33333,-2.3094) to [bend right=40] (-2,-3.4641);
\draw[dashed] (-1.33333,-2.3094) to [bend left=40] (-2,-3.4641);
\draw (4,0) to [bend right = 40] (2.66667,0);
\draw[dashed] (4,0) to [bend left = 40] (2.66667,0);
\draw (-12,0) ellipse (3 and 1.5);
\begin{scope}
\clip (-12,-1.8) ellipse (3 and 2.5);
\draw (-12,2.2) ellipse (3 and 2.5);
\end{scope}
\begin{scope}
\clip (-12,2.2) ellipse (3 and 2.5);
\draw (-12,-2.2) ellipse (3 and 2.5);
\end{scope}
\filldraw[red] (-14.5,0) circle (3pt) node [above] {\footnotesize{$2$}};
\filldraw[red] (-9.5,0) circle (3pt) node [above] {\footnotesize{$2$}};
\filldraw[red] (-12,-1.25) circle (3pt) node [above] {\footnotesize{$2$}};
\filldraw[red] (0,0) circle (3pt);
\draw[->, shorten <=2cm, shorten >=2cm] (-12,0) -- (0,0);
\node at (2,2.5) {\footnotesize{$X_3$}};
\node at (-12, 2.5) {\footnotesize{$D_3$}};
\end{tikzpicture}
\caption{The immersion of $D_3$ in $X_3$}\label{fig:EmbedD3}
\end{figure}

Since $p \notin \Gam_1, \Gam_2$, the cone point $z_2 \in D_0$ of order two, whose image in $X_0$ has orbifold weight four, has nine preimages in $X_3$ that are all smooth points for the orbifold structure. Since $z_2$ also has exactly nine preimages in $D_3$, it follows that $D_3 \looparrowright X_3$ is an embedding on these points and thus $D_3$ meets $D_3^\prime$ at each of these points, and each intersection is transversal since this is the case in the universal cover. Since $q \in \Gam_1$ but $q \notin \Gam_2$, the cone point $z_3$ of order three has $\frac{36}{3}=12$ lifts to $X_3$ and $\frac{18}{3}=6$ lifts to $D_3$. The action of $(q, \Id)$ at $\wt{z}_3$ in Figure \ref{fig:Diag0} identifies the three embedded copies of $\bbH^2$, which means that $D_3$ does not meet $D_3^\prime$ at these points, so the twelve preimages are distributed as six smooth points on each of $D_3$ and $D_3^\prime$. Thus $D_3$ embeds in $X_3$ away from its three cone points of order two, as claimed.

The cone point $z_{12}$ on $D_0$ of order twelve has preimage the three cone points of order two on $D_3$. The image of $z_{12}$ on $X_0$ has orbifold weight $12^2=144$, and it lifts to a single orbifold point of weight four with local group $(\bbZ / 2)^2$ on $X_3$ associated with the conjugacy class of the subgroup of $\Gam_1 \times \Gam_2$ generated by $(r^6, \Id)$ and $(\Id, r^6)$. Indeed, this is the point $(w_1, w_2)$ on $X_3$ where $w_j$ is the unique cone point of order two on $\Gam_j \bs \bbH^2$. The cone points of order two on $D_3$ then must all map to this point, hence the image of $D_3$ passes through this point with multiplicity three. This completely determines the image of $D_3$ in $X_3$, which fails to be an embedding at precisely the three cone points of order two, which are all identified. Thus the image of $D_3$ is as depicted in Figure \ref{fig:EmbedD3}, proving the lemma.
\end{pf}

\begin{lem}\label{lem:SIstep3}
Ramification of the orbifold cover $Y \to X_3$ is:
\begin{itemize}

\item degree $n$ ramification above vertical curves of the form $\{x_n\} \times \Gam_2 \bs \bbH^2$, where $x_n \in \Gam_1 \bs \bbH^2$ is a cone point of order $n$;

\item degree two ramification over $\Gam_1 \bs \bbH^2 \times \{w_2\}$, where $w_2$ is the cone point of order two on $\Gam_2 \bs \bbH^2$.

\end{itemize}
If $p_j$ denotes projection of $D_3$ onto $\Gam_j \bs \bbH^2$ and the image $(p_1(z), p_2(z))$ of a point $z \in D_3$ on $X_3$ intersects the ramification locus, then $p_j(z)$ is a cone point of $\Gam_j \bs \bbH^2$ for at least one $j \in \{1,2\}$.
\end{lem}

\begin{pf}
The last statement is an immediate consequence of the first and Lemma \ref{lem:SIstep2}. For the first statement, note that ramification comes from elements of the form $(t, \Id)$ or $(\Id, t)$ not contained in the subgroup $\Lam$ of $\Gam_1 \times \Gam_2$ associated with $Y$, where $t$ is finite order. Indeed, elements of $\Gam_1 \times \Gam_2$ with nontrivial entry in each coordinate have an isolated fixed point on $\bbH^2 \times \bbH^2$, so only elements with one element trivial can act as a reflection through a vertical or horizontal hyperbolic planes, which is precisely what creates ramification in $Y \to X_3$. Since $\Lam$ contains none of the elements of $\Gam_1 \times \Gam_2$ with this form by definition, as $\rho_1$ and $\rho_2$ are injective on torsion, the lemma follows.
\end{pf}

\begin{pf}[Proof of Proposition \ref{prop:SelfIntersection}]
Set $D_4 = \Gam_4 \bs \bbH^2$ and let $\wh{E}$ denote its image in $Y$. The map $D_4 \to D_3$ has degree four, and $Y \to X_3$ has degree $12$ and is unramified over $D_3$ by Lemma \ref{lem:SIstep3}. Thus $D_3$ has three preimages in $Y$.

Since $X_3$ is a product and $D_3$ embeds away from its cone points of order two by Lemma \ref{lem:SIstep2}, the image of $D_3$ in $X_3$ meets each vertical curve $\{z_3\} \times \Gam_2 \bs \bbH^2$, over which the map $Y \to X_3$ has degree three ramification by Lemma \ref{lem:SIstep3}, in a single point that has four preimages in $Y$. Since $D_4 \to D_3$ is of degree four and unramified above the relevant points, $\wh{E}$ must pass through all four preimages. Thus $\wh{E}$ is smooth through each lift.

Again by Lemma \ref{lem:SIstep2}, the image of $D_3$ meets both $\{x_2\} \times \Gam_2 \bs \bbH^2$ and $\Gam_1 \bs \bbH^2 \times \{w_2\}$ precisely in the point $(x_2, w_2)$ where the immersion of $D_3$ fails to be an embedding. This point on $X_3$ has preimage containing the two $\mathrm{A}_1$ singularities of $Y$. Locally around an $\mathrm{A}_1$ singularity the cover is the map from a neighborhood $V$ of an $\mathrm{A}_1$ singularity to a neighborhood $U$ of the origin in $\bbC^2$ given by
 \[
 U \lra V \lra U
 \]
 where the composition $U \to U$ is the action of $(\bbZ/2)^2$ by reflections of order two in each coordinate and the map $U \to V$ is the quotient by the diagonal subgroup generated by $(-1,-1)$. The map $V \to U$ is two-to-one on lines through the singular point. It follows that the preimage of $D_3$ in $Y$ passes through each $\mathrm{A}_1$ singularity with multiplicity six.

On the other hand, each of the twelve cone points on $D_4$ must map to one of the $\mathrm{A}_1$ singularities, since these are the only orbifold points on $Y$ for its $\bbH^2 \times \bbH^2$ orbifold structure, and orbifold points of $D_4$ must map to orbifold points of $Y$. It follows that $\wh{E}$ must pass through the $\mathrm{A}_1$ singularities with total multiplicity twelve, and therefore it passes through each singularity with multiplicity six. This proves that $\wh{E}$ is embedded in $Y$ away from the $\mathrm{A}_1$ singularities and passes through each singular point with multiplicity six.

Now consider the proper transform $E$ of $\wh{E}$ to $Z$ under the minimal resolution $\varphi : Z \to Y$, and let $F_1, F_2$ be the $(-2)$ curves resolving the singularities. Since $E$ meets each singular point with multiplicity six,
\[
E \cdot F_j = 6
\]
for $j = 1,2$. Write $E$ as $\varphi^* \wh{E} + \al_1 F_1 + \al_2 F_2$ for $\al_j \in \bbQ$. Then $K_Z \cdot F_j = 0$, so
\begin{align*}
K_Z \cdot E &= \varphi^* K_Y \cdot \varphi^* \wh{E} \\
&= K_Y \cdot \wh{E} \\
&= 2|e^{\mathrm{orb}}(D_4)|
\end{align*}
by relative Hirzebruch proportionality \cite[\S 4]{HirzebruchRel}, where $e^{\mathrm{orb}}$ denotes the orbifold Euler characteristic. Since
\[
e^{\mathrm{orb}}(D_4) = 0 - 12 \times \frac{1}{2} = -6
\]
it follows that $K_Z \cdot E = 12$. Since $E$ has genus one, adjunction implies that $E^2 = -12$, which completes the proof of the proposition.
\end{pf}

\subsection{Proof that $(Z, E)$ is a ball quotient pair}\label{ssec:BallQ}

Let $(Z, E)$ be the pair constructed in \S\ref{ssec:FindElliptic}. First, there are the simple calculations
\begin{align*}
c_1^2(Z, E) &= (K_Z + E)^2 = K_Z^2 - E^2 = 18 \\
c_2(Z, E) &= c_2(Z) = 6
\end{align*}
that give $c_1^2(Z, E) = 3 c_2(Z, E)$. As described in \S \ref{sec:Background}, to prove that $(Z, E)$ is a ball quotient pair, it suffices to prove that $K_Z + E$ is nef \cite[Def.\ 1.4.1]{Lazarsfeld} and big \cite[Def.\ 2.2.1]{Lazarsfeld}, and that every $(-2)$ curve of $X$ meets $E$ nontrivially.

Consider an irreducible curve $C$ on $Z$. If $C$ is not $E$ or one of the $(-2)$ curves $F_j$, then $C \cdot E \ge 0$ and $\varphi(C)$ is an irreducible curve on $Y$. Now note that $K_Y$ is ample, since $Y$ is a compact quotient of $\bbH^2 \times \bbH^2$ by a group action with only isolated fixed points. Indeed, since there is no branching divisor, the orbifold canonical class for $Y$ is the same as the canonical class \cite[Prop.\ 4.4.15]{BoyerGalicki}. Then
\[
(K_Z + E) \cdot C = K_Y \cdot \varphi(C) + \wh{E} \cdot C > 0
\]
since $K_Y$ is ample.

From this, we obtain that
\begin{align*}
(K_Z + E) \cdot E &= 0 \\
(K_Z + E) \cdot F_j &= 6
\end{align*}
so $(K_Z + E) \cdot C \ge 0$ for all irreducible curves $C$ on $Z$, and therefore $K_Z + E$ is nef. Then $(K_Z + E)^2 > 0$, so it is big by \cite[Thm.\ 2.2.16]{Lazarsfeld}. Moreover, note that ampleness of $K_Y$ implies that there is no $(-2)$ curve on $Y$ disjoint from the singular points, since such a curve $C$ would have $K_Y \cdot C = 0$ by adjunction, contradicting ampleness of $K_Y$. Therefore there is no $(-2)$ curve on $Z$ that does not meet $E$. This combined with Chern--Gauss--Bonnet completes the proof of the following result:

\begin{thm}\label{thm:FirstExDone}
The pair $(Z, E)$ constructed in this section is a ball quotient pair. Specifically, $Z \ssm E$ is a smooth one-cusped ball quotient of volume $16 \pi^2$.
\end{thm}

\section{The proof of Theorem \ref{thm:MainTech}}\label{sec:Bootstrap}

Recall that Theorem \ref{thm:Main} is an immediate consequence of Theorem \ref{thm:MainTech} and Chern--Gauss--Bonnet. This section proves Theorem \ref{thm:MainTech}. To start, the pair $(Z_1, E_1)$ for Theorem \ref{thm:MainTech} is the surface $Z$ from \S \ref{sec:Construction} and elliptic curve $E$ on $Z$, so it suffices to construct $(Z_d, E_d)$ for all odd $d > 1$.

\begin{lem}\label{lem:CoverReduction}
To prove Theorem \ref{thm:MainTech}, it suffices to prove that there is an \'etale cover $Z_d \to Z_1$ of every odd degree $d$ so that $E_1$ has irreducible preimage $E_d$ in $Z_d$.
\end{lem}

\begin{pf}
Note that $E_d \to E_1$ is \'etale, so $E_d$ has genus one. Given the hypotheses of the lemma, there is an \'etale covering
\[
Z_d \ssm E_d \lra Z_1 \ssm E_1
\]
of degree $d$. Since $Z_1 \ssm E_1$ is a ball quotient, so is $Z_d \ssm E_d$. Since $E_d$ is a single genus one curve, $Z_d \ssm E_d$ is a one-cusped ball quotient and hence $(Z_d, E_d)$ is a ball quotient pair. Moreover, characteristic classes multiply in covers, so $c_1^2(Z_d) = c_2(Z_d) = 6 d$. Self-intersection also multiplies, so $E_d$ has self-intersection $-12d$. Thus the surfaces $(Z_d, E_d)$ satisfy the hypotheses of Theorem \ref{thm:MainTech} and hence their existence would prove that result.
\end{pf}

\begin{prop}\label{prop:Covers!}
For each odd $d \ge 1$ there is an \'etale cover $Z_d \to Z_1$ of every odd degree $d$ so that $E_1$ has irreducible preimage $E_d$ in $Z_d$.
\end{prop}

\begin{pf}
Recall from Proposition \ref{prop:DescribeQuo} that the Albanese map of $Z_1$ is given by projection onto the smooth curve $D_2 = \Gam_2 \bs \bbH^2$ of genus one. It follows from standard facts about the Albanese map that $\pi_1(Z_1)^{\mathrm{ab}}$ modulo torsion is isomorphic to $\bbZ^2$ induced by the Albanese map (e.g., see \cite[\S I.13]{BPV}). By construction, the map $E_1 \to D_2$ is an \'etale cover of degree $12$ (e.g., recall Figure \ref{fig:Covers2}), hence the image of $\pi_1(E_1)$ in $\bbZ^2$ has index $12$. Thus for each odd $d \ge 1$ there is a finite abelian quotient $A_d$ of $\bbZ^2$ with order $d$ so that the induced homomorphism $\mu_d : \pi_1(Z_1) \to A_d$ has restriction to $\pi_1(E_1)$ that is surjective. If $Z_d$ is the covering of $Z_1$ associated with $\mu_d$, then $E_1$ has exactly one preimage in $Z_d$, which then must be an irreducible curve of genus one. This proves the proposition.
\end{pf}

This completes the proof of the main result of this paper.

\begin{pf}[Proof of Theorem \ref{thm:MainTech}]
Take the surfaces $(Z_d, E_d)$ provided by Proposition \ref{prop:Covers!} and apply Lemma \ref{lem:CoverReduction}.
\end{pf}



\end{document}